\documentclass[a4paper,11pt]{article}
\usepackage[utf8]{inputenc}
\usepackage[T1]{fontenc}
\usepackage[portuguese]{babel}

\usepackage[]{graphicx}
\usepackage{geometry}
\usepackage{verbatim}
\usepackage{latexsym,amsfonts,amssymb}
\usepackage[normalem]{ulem}
\usepackage{color}
\usepackage{subfigure}
\usepackage{varioref,float,amsmath,amsthm}
\usepackage{mathtools}
\usepackage{epstopdf}

\begin{document}

\begin{center}
.\\\vspace{7cm}
\noindent\textbf{\LARGE An overview of the fractional-order gradient descent method and its applications}\\ 

\vspace{2cm}
{\small\noindent{ Higor V. M. Ferreira$^{1}$, Camila A. Tavares$^{2}$, Nelson H. T. Lemes$^{1,*}$ 
and Jos\'e P. C. dos Santos$^{3}$}} \\
\vspace{2cm}
\end{center}

\noindent
\begin{center}
{\small 
$^{1}$ Laboratory of Mathematical-Chemistry, Institute of Chemistry, Universidade Federal de Alfenas (UNIFAL), Alfenas, MG, Brazil  \\
$^{2}$  
Institute of Chemistry and Center for Computational Engineering \& Science,
Universidade Estadual de Campinas (UNICAMP), Campinas, SP, Brazil\\
$^{3}$ Institute of Exact Science, Universidade Federal de Alfenas (UNIFAL), Alfenas, MG, Brazil  \\
}
\end{center}
\vspace{\stretch{1}}
\noindent{$^*$nelson.lemes@unifal-mg.edu.br}

\newpage
\subsubsection*{Abstract}
Recent studies have shown that fractional calculus is an effective alternative mathematical tool in various scientific fields. However, some investigations indicate that results established in differential and integral calculus do not necessarily hold true in fractional calculus. In this work we will compare  various methods presented in the literature to improve the Gradient Descent Method, in terms of convergence of the method, convergence to the extreme point, and convergence rate. { In general, these methods that generalize the gradient descent algorithm by replacing the gradient with a fractional-order operator are  inefficient in achieving convergence to the extremum point of the objective function.} 
To avoid these difficulties, we proposed to choose the Fractional Continuous Time algorithm to generalize the gradient method. In this approach, the convergence of the method to the extreme point of the function { is guaranteed by introducing the fractional order in the time derivative, rather than in of the gradient.}
In this case, { the issue of finding the extreme point is resolved, while the issue} of stability at the equilibrium point remains.
{ Fractional Continuous Time method converges to extreme point of cost function 
when  fractional-order  is between 0 and 1. 
The simulations shown in this work suggests that a similar result can be found when 
$1 \leq \alpha \leq 2$.}
{ This paper highlights the main advantages and disadvantages of generalizations of the gradient method using fractional derivatives, aiming to optimize convergence in complex problems. Some chemical problems, with {$n=11$ and $24$} optimization parameters, are employed as means of evaluating the efficacy of the propose algorithms. In general, previous studies are restricted to mathematical
questions and {simple illustrative examples}.} \\

\noindent
{ Keywords:} Fractional calculus; Descent gradient method; Chemical optimization problem

\clearpage

\section{Introduction}

Fractional calculus is a generalization of classical calculus, in which preserves many of its basic properties, but, also offers new features for research, such as the memory effect \cite{3j}.
Because of this, it has been increasingly used in various applications \cite{Nel1,Nel2,Nel3,Nel4}. Despite the growing use of this tool in science, many fundamental questions still remain open. For example, there are many different definitions of fractional derivative operators, which do not coincide in general. Each one of them 
try to preserve different properties of the classical derivative. The Riemann-Liouville, Caputo e Gr\"unwald-Letnikov are the three most commonly used definitions for fractional differentiation $D_u^\alpha f(u)$ \cite{3j}. In physical research, in general, the emphasis is given to the Caputo definition, because it provides initial conditions with clear physical interpretation for the differential equations of the fractional order \cite{3j}.

It is well-known that integer-order derivatives and integrals have clear physical interpretation and are used for describing different concepts in science.  However, for the fractional derivative, the situation is complicated  and no physical or geometric interpretation has been fully used until now. 
Accordingly, we cannot  
assume 
{ that our experience with integer-order calculus will apply in fractional calculus.}
For example, the first-order derivative equal to zero means that we have found the critical point. In general, the same statement cannot be made about fractional derivative order 
leading
to unexpected consequences when compared with the integer-order derivative \cite{8P,9P}. In this work,
we will explore the use of fractional-derivative order in the gradient descent method (GDM). 
This is the main challenge raised against this line of research, { and a solution is presented here.}

GDM is one of the simplest and most commonly used methods to solve problems in science and engineering involving finding a set of parameters ${ u}$ that minimizes an objective function $f({ u})$.
The original gradient descent algorithm is attributed to Cauchy, who first suggested it in 1847 \cite{H1}. 
However, there are a number of
variations, some of which are very recent \cite{Ruber,Ruber2}. 
An overview of the advantages and disadvantages of GDM can be found in the literature.  See, for example, Ruder et al. (2016) \cite{Ruber}.
Although the method has been successful in many problems, time convergence is a property that can always be improved, because in general, the gradient descent method is slow to converge close to the extreme point \cite{Ruber}. In this perspective, fractional
calculus {has been  little explored as an alternative} to improve the ordinary gradient descent method (also known as the steepest descent method) \cite{Rede1,Rede2,Rede3}, due to the non-local characteristics of the fractional derivative, which  is an important contrast to the derivative of integer order \cite{3j}. 

The GDM is a way to minimize an objective function $f(u)$, over a set of parameters $u$,
by updating these parameters in the opposite direction of the derivative of the integer order of the $f$,  with respect to $u$, at $u$ (this is the gradient of the $f$ in one-dimensional space). Therefore, the gradient descent method starts with an initial value of $u$ and constructs a sequence of values of $u$ such that $f(u_{i})<f(u_{i-1})$ at every interaction $i$. The expected result is that the sequence of $u$ converges to an extreme point $u^*$ of the function $f(u)$, such that $df/du=0$ at $u=u^*$. 
In practice, a stopping criterion is used to decide when to stop the interaction process.
Therefore, in the gradient descent method, the equilibrium point $u^\#$ (such as $du/dt=0$) agrees with the extreme point $u^*$ (such as $df/du=0$), since $du/dt=-\rho df/du$ in the
gradient descent  method. Thus, the initial optimization problem is reduced to the problem of finding an  equilibrium point. This is the central point of the whole procedure.

{ The first attempt to generalise GDM to fractional order is given in reference \cite{4P}. The authors modify the 
gradient descent method using the Riemann-Liouville derivative of fractional order  instead of the gradient of $f$ (i.e. the first-order derivative of $f$). However, in this case, it cannot  be guaranteed that the convergence point {(when $du/dt\rightarrow 0$ and} $u(t>t^*)=$constant) provides an extreme point defined by $df/du=0$, whereas $_0D_u^\alpha f=0$  does not define an extreme point of the objective function $f(u)$.
A similar idea can be found in reference \cite{5P}, where the Caputo derivative of fractional order was used instead of the gradient of the $f$.
Once again, it cannot be guaranteed that the final result obtained by this method
is the same that provides the extreme point of the objective function $f(u)$.
In other words, for both derivative operators, when $_0^*D_u^\alpha f=0$ at 
$u=u^\#$ does not imply that $u^\#$ is an extreme point of the $f(u)$. Although
$u^\#$ represents a equilibrium point, since $du/dt=-\lambda _0^*D_u^\alpha f=0$ at 
$u=u^\#$. This holds for any definition chosen for $D_u^\alpha f(u)$ among the classical fractional derivatives. 
A more detailed discussion of this point is presented in the methodology section.}

To circumvent this difficulty, in reference \cite{6P} a new modification was proposed. The authors used the Riemann-Liouville and Caputo derivatives of fractional order  with a fixed memory length. Now it can be shown that the proposed method converges to the proximity of the extreme point when memory length is small.
{Now $_LD_u^\alpha f=0$ at $u=u^\#$ implies that $u^\#$ is a point close to the extreme point $u^*$ of the $f(u)$.}
{However,
when memory length is small, the non-local characteristics of the fractional derivative are lost. Therefore, the proposed model returns, in a certain way, to the
original method.}  It is noteworthy that the solution found by the previous  model, reference \cite{6P}, was achieved in less time when compared to the original method, which is a great result.
Although the fractional gradient method, with a fixed memory length, gained competitive advantages over the integer-order gradient method, the convergence to the extreme point {is only approximately guaranteed} \cite{7P}. Therefore,  this question should be further explored. In reference \cite{6P}, the authors said that  ``research of the fractional-order gradient method (FDGM) is still in its infancy and deserves further investigation".

{This paper presents a different approach and its applications of incorporating the fractional derivative operator into the iterative algorithm of the gradient descent method.
The approach proposed in this paper differs from previously attempted strategies, in \cite{4P,5P,6P}.
The basic idea of this paper is to replace the first-order derivative in time  by a non-integer derivative operator, 
$D_t^\alpha u(t)=-\rho df(u)/du=F(u)$.  
This procedure is called of the Fractional Continuous Time Method (FCTM) \cite{ine12}.  
Here, the parameters are still updated in the opposite direction of the gradient of $f$, with respect to $u$, evaluated at $u$.
The values of $u$ where $F(u) = 0$ are called equilibrium points of the dynamic system. In this case, $F(u) = 0$ implies $\frac{df(u)}{du} = 0$. 
Therefore, the equilibrium points of the FDE are the extreme points of the function $f(u)$. 
The reference \cite{O2} provides a proof 
that the proposed method converges to the exact value of the extreme point of the $f(u)$ when $0<\alpha<1$. 
In this paper, the simulations suggest that a similar result can be found when 
$1 \leq \alpha \leq 2$.
Besides this, good approximations have been achieved in less time when compared to the usual GDM. Therefore, {the proposed method} keeps the advantages of the previous method \cite{4P,5P,6P} and
corrects the failures of the previous approach, ensuring convergence to the exact value of the extreme point. The disadvantage is that the proposed model converges more slowly than the 
Cauchy method after several iterations.}

{Therefore,
the main goal of this paper is  to analyze the efficiency of FCTM in chemical  optimization problems.
The majority of  applications in literature \cite{4P,5P,6P,7P,ine12} are restricted to the optimization of problems with a limited number of variables, typically $n=1$ or $2$. The present study will examine both, linear and non-linear optimization problems, each involving {more than eleven} variables.
The results obtained in two cases, by our model, will be compared with the common integer order model.
The first example
deals with the linear system
in which ${ K}$ is a
Vandermonde matrix, with $Dim({ K})=11\times 11$ \cite{Nel1,Nel2}.
Finally, {the proposed model}  will be  used to
find the  minimum energy configuration of $N$ point charges on a surface of the unit sphere. This problem originated with Thomson's ``plum pudding" model of the atom \cite{1T}. After considerably more structural information, this idea was  extended by Gillespie and Nyholm in the Valence-Shell-Electron-Pair-Repulsion (VSEPR) theory \cite{2T}. The Thomson is a classical, well-known and
yet very important problem within the Physical-Chemistry and optimization  methods research \cite{3T,4T,5T,6T,7T}.

{The rest of the paper is organized as follows: 
In the next Section 2, we will begin by introducing our generalized model and the necessary foundations of
fractional calculus will be presented.
The prototype problems for the examination of the method proposed will be presented in Section 3, together with a discussion
of the efficiency of the strategy proposed here.
Some conclusions will be
included in Section 4.}

\section{Methodology}


\subsection{Fractional calculus background}

Since 1695, many possible definitions have been proposed for the fractional operator, among them three are widely used: the Gr\"unwald$-$Letnikov, the Riemann$-$Liouville, and the Caputo. 
These definitions are different in the domain of the function $f(u)$.
The Gr\"unwald-Letnikov $\alpha$th-order fractional derivative of a function $f(u)\in \mathbb{R}^n$ with respect to $u\in(0,u_{max}]$ is given by  a generalization of the formula for the  $n$th
{derivative of $f(u)$,} with $n$ = 1, 2, 3 ... \cite{3j},
\begin{equation}
{
_a\mathcal{D}_u^{(\alpha)} f(u)\coloneqq
\lim_{h\rightarrow 0^+} \frac{1}{h^\alpha}
\sum_{k=0}^{N} (-1)^k \left(
\begin{array}{c}
\alpha\\
k
\end{array}
\right) f(u-kh)
}
\label{eqdglx}
\end{equation}
where $n$ was replaced by fractional-order $\alpha$ ($\alpha \in \mathbb{R}^+$) and the binomial coefficient was written as, 
\begin{equation}
\left(
\begin{array}{c}
\alpha\\
k
\end{array}
\right)=C_{\alpha,k}=
\frac{\Gamma(\alpha+1)}{\Gamma(k+1)\Gamma(\alpha-k+1)}
\label{bincoeffx}
\end{equation}
{In this case, Euler's gamma function is the choice for the interpolation of the factorial function.} 
{Therefore, the equation (\ref{eqdglx}) is well defined for all  positive real number $\alpha$. }  
The total terms in the sum on the Equation (\ref{eqdglx}), $N=\lceil (t-a)/h\rceil$, is the smallest integer such that $N > t/h$. The sum shown in Equation (\ref{eqdglx}) converges absolutely and uniformly for all $\alpha > 0$ and for every bounded function $f(u)$ \cite{3j}.

In order to
present Riemann-Liouville fractional derivative, we will start
with the generalization of the integral of
non-integer order, so if we
integrate $f(u)$
{with respect to $u$ by}
$m$ times we find the Cauchy formula for repeated integration
\cite{3j},
Now, using the generalization of factorial by Gamma function,
we define the fractional integral of Riemann$-$Liouville,
\begin{equation}
J^{\alpha}f(u)
\coloneqq \frac{1}{\Gamma(\alpha)}\int_a^u (u-s)^{\alpha-1} f(s) ds
\label{eqint2}
\end{equation}
where now $m\geq 0$
and $-\infty\leq a\leq u\leq +\infty$.

The definition of fractional derivative of Riemann-Liouville can be obtained using the
fundamental theorem of calculus,
in which
$D^rJ^r f=f$.
{In order to suggest a fractional
derivative definition, consider that $r=q-\alpha$, $D^{q-\alpha}J^{q-\alpha} f=f$.
Now, apply the $D^\alpha$ operator to both sides of the equation,
$D^\alpha D^{q-\alpha}J^{q-\alpha} f=D^\alpha f$.}
{If the composition rule was valid for non-integer exponents, we would have an expression for $D^\alpha f$, 
which suggests defining the fractional derivative as}
\begin{equation}
{_aD_u^\alpha f(u)\coloneqq D^nJ^{n-\alpha }f=\frac{1}{\Gamma(n-\alpha)}\frac{d^n}{du^n}\int_a^u\frac{f(s)ds}{(u-s)^{\alpha -n+1}}}
\label{eqdif1}
\end{equation}
{where
${_aD_u^0 f(u)=f(u)}$, $n\in \mathbb{Z}^+$, $n-1\leq \alpha\leq n$ and $u>a>0$.}
Since $J^{n-\alpha }$ is a well-defined operator when $n-\alpha$ is a noninteger number, 
the Equation (\ref{eqdif1}) represents a well defined
fractional derivative operator, called the Riemann-Liouville fractional derivative. 

{Another attempt to define the fractional derivative was made by Caputo in 1967 \cite{2j,3j,4j},
who suggested changing the order of the operators $D^n$ and $J^{n-\alpha }$ in {Equation (\ref{eqdif1})}.}
Thus, Caputo fractional derivative
is defined
as
\begin{equation}
{^*_aD_u^\alpha f(u)\coloneqq J^{n-\alpha }D^nf=\frac{1}{\Gamma(n-\alpha)}\int_a^u\frac{\frac{d^nf(s)}{ds^n}ds}{(u-s)^{\alpha -n+1}}}
\label{eqdif}
\end{equation}
{in which, ${^*_aD_u^0 f(u)=f(u)}$, $n-1\leq \alpha\leq n$ and $u>a>0$.}
 In the same way that $JDf\neq DJf$, we found that ${^*_0D_u^\alpha f\neq {_0D}_u^\alpha f}$.
From the three definitions above, the Caputo definition is the
most often used in Physical-Chemistry applications. The main advantage of the Caputo
derivative is that it only requires initial conditions given in terms of
integer-order derivatives, representing well-understood features
of physical situations and thus making it more applicable to real
world problems \cite{2j,3j,4j}.

The fractional operator is clearly non-local, since the fractional derivative depends on the lower boundary of the integral.
This is a contrast to the derivative of integer order,  clearly a local operator.
Another
properties of Caputo fractional derivative are \cite{2j}:
a) the Caputo fractional derivative of order zero returns the input function; b) the
Caputo fractional derivative of integer order $n$ gives the same result as the usual differentiation of order $n$;
c) {the Caputo fractional derivative of a constant function is equal to zero }
and d) the Caputo fractional derivative is an ill-posed  operator, this means that small errors in input data may yield large errors in the output result.
{ The classical fractional derivatives are those of Caputo, Riemann-Liouville, and Gr\"unwald-Letnikov. These classical fractional derivatives satisfy the criteria proposed by Ortigueira and Tereiro Machado in 2015 
\cite{crit}.}
For a complete review of the Fractional Calculus see, for example, the book of Podlubny \cite{3j}.

\subsection{Fractional descent gradient method}

\subsubsection{Original method}

Before the introduction of fractional gradient descent (FGDM)
method, it is worthwhile to retrospect the gradient descent (GDM)
method within the integer-order framework. See, for example, Lemarechal (2012) \cite{H1}. 
Often, this model has been used to find the solution ${u}^*$ that minimizes an objective function $f(u)$. It is well-known that iterative step of the conventional gradient method, in one dimensional problem, is given by 
\begin{equation}
u_{k+1}=u_{k}- \left.\omega\hspace{0.1cm} \frac{d}{du}f(u)\right|_{u=u_k}
\label{gradeq}
\end{equation}
{in which $\omega$ is the step size and $k$ is the number of iterations}. More details can be found 
in reference \cite{Ruber,Ruber2}. Briefly, the
method can be described as: the new variable at step $k$ is updated using
the  first-order derivative of $f$
in relation a $u$. 
Using appropriated $\omega$, 
the result of the gradient descent algorithm is a monotonic sequence
$f(u_0)>f(u_1)>...>0$. By this method is expected that the sequence $u_k$ converges to a local minimum $u^*$ of the function $f(u)$.

The GDM can be seen, by the explicit Euler method, as a discretized form of the following ordinary differential equation,
\begin{equation}
\frac{du}{dt}=- \lambda\frac{d}{du}f(u)=F(u)
\label{diferentialgradeq}
\end{equation}
in which $t$ is the iteration time and $\omega=\lambda dt$.  Now, the method above is called Continuous Gradient Method (CGM). 
In this model,
the variable $u$ propagates until it reaches an equilibrium point $u^\#$, when $F(u)= 0$. {By considering the above equation, it follows that the equilibrium point 
$u^\#$ represents an extreme point $u^*$ of the $f(u)$, when $\frac{d}{du}f(u)$ is also equal to zero.}

{Some changes have been made in Equation (\ref{gradeq}) to include fractional derivative order on the direction in which the parameter $u$ is updated \cite{Ruber,Ruber2}. When proposing new modification of the GDM, some 
questions must be answered: a)If we start with an initial value of $u_0=u(0)$ and build the sequence of values $u_k$, according to the above rule, do we have $f(u_k) < f (u_{k-1})$ for every iteration $k$? b) Is the sequence $u_k$, which is obtained by the above rule, convergent?
c) The sequence $u_k$ obtained by the above rule converges to the extreme point $u^*$ of the function $f(u)$ such that $df/du = 0$ at $u = u^*$? and d) The sequence, $u_k$, obtained by the above rule converges faster than the sequence obtained by the original Cauchy method? 
Some of them will be detailed in the following subsection, together with the method 
that will be discussed here.} 

\subsubsection{Previous work}

In references \cite{4P,5P}, the authors suggest improving the gradient method by replacing
the first-order derivative with respect to $u$ by fractional-order derivative. 
Therefore, the
Equation (\ref{gradeq}) was rewritten as
\begin{equation}
u_{k+1}=u_{k}- \omega\hspace{0.1cm} [_0D^\alpha_u f(u)](u_k) 
\label{liteq}
\end{equation}
in which the Riemann-Liouville definition of fractional derivative was used in reference \cite{4P} and Caputo operator in reference \cite{5P}.
For integer-order algorithms, it 
is well-known that
{the sequence $u_k$ converges to an equilibrium point $u^\#=u^*$, where $u^*$ is a extreme point of $f$.}
{When we change integer derivative  by fractional derivative order, the equilibrium points are changed, i.e.
$_0D_u^\alpha f(u)|_{u=u^*}\neq 0$. Now the value of $u$ that $_0D_u^\alpha f(u)= 0$ is
different to $u^*$.}
Therefore, the solution of
integer-order differential equation (\ref{gradeq}) is different from the one obtained by fractional differential equation (\ref{liteq}).  This is 
the main disadvantage of the {algorithm \cite{4P,5P}}. 

For example, considering the function $f(u)=(u-c)^2$, in which extreme point is $u^*=c$, 
can be easily determined that $_0^*D_u^\alpha f(u)=0$ when $u^\#=c(2-\alpha)$. 
In this case, $u^\#$ is close to $c$ when $\alpha$ is close to $1$.
As a result, the value of $u$, such as $_0^*D_u^\alpha f(u)=0$,  is not the extreme point of the function $f(u)$.
{ Figure (\ref{figure1}) shows \( u(t) \) at each iteration obtained by the FGDM method, using the Caputo derivative with \( \alpha = 0.9 \). In this case, it can be observed that \( u(t) \) does not converge to the extreme point of the function \( f(u) \).}

Note that this result would be different if the Riemann-Liouville definition  of fractional derivative is used. By using the Riemann-Liouville definition, we arrive at the following equation
\begin{equation}
\frac{\Gamma(3)}{\Gamma(3-\alpha)}u^2-2c\frac{\Gamma(2)}{\Gamma(2-\alpha)}u+c^2 \frac{\Gamma(1)}{\Gamma(1-\alpha)}=0
\end{equation}
The above equation has two roots. Consequently, there are two values of $u$ such as 
$_0D_u^\alpha f(u)=0$. 
{Therefore, by replacing $d/dt$ with ${_aD}_t^\alpha$, {in equation (\ref{liteq})}, the result obtained depends on the chosen fractional derivative operator,  Caputo or Riemann-Liouville.} 

{Based on the fractional integral with the parametric lower limit of integration, $a$, 
the following are the definitions of}
Riemann-Liouville and Caputo.
In addition to this, the value of $u$, such that $_a^*D_u^\alpha f(u)=0$, depends on the value of $a$. 
Consequently, only the function on the interval $[a, t]$ is taken into account by the fractional derivative operator.
Considering the function $f(u)=(u-c)^2$, the value of $u$, such that $_a^*D_u^\alpha f(u)=0$,  is given by
\begin{equation}
u=a+c(2-\alpha)
\end{equation}
However, this issue were omitted in the discussion in the references \cite{4P,5P}. See references \cite{8P,9P} for a detailed discussion of extreme points with
fractional-order derivatives.  
Nevertheless, these examples raise the need 
for more research in the method proposed. 

{Despite the  fact that the  method proposed in reference \cite{4P,5P} cannot guarantee convergence to an extremum, the method predicts sufficiently close results in less time than the integer-order model \cite{4P,5P}. This is an important advantage of this method  compared to the traditional approach, which motivates further 
 interest in studying the generalized algorithm. Thus, the generalized algorithm
can be used to generate better initial conditions, which will be utilized in algorithms where convergence to the extremum is guaranteed. }

To circumvent the difficulties mentioned above, the paper \cite{6P}
uses fractional-order derivative with limited memory length. In that case, the
generalized model was rewritten as
\begin{equation}
u_{k+1}=u_{k}- \omega\hspace{.1cm} [_aD_u^\alpha f(u)](u_k) 
\label{liteq2}
\end{equation}
where $a$ is  
lower limit of integration on fractional derivative operator. 
However, when the memory length is short, the fractional-order derivative becomes nearly integer-order derivative \cite{10P}. 
{With the suggested modification, in \cite{6P}, the nonlocal operator becomes a local operator.} To a certain extent, the proposed modification returns to the original model. 
Therefore, it is expected that the value of  $u$ for which  
$_aD_u^\alpha f(u)=0$ stays close to the value of $u$ for which $df/du=0$. Even in this case, the
convergence to the extreme point cannot be guaranteed \cite{7P}. { Figure \ref{figure1} also shows the result obtained by this method. As we can observe, the solution for this method converges near the extreme point.}

Considering that the fractional 
derivative of the Caputo type  is given by the following Taylor series expansion \cite{Taylor,Taylor2}
\begin{equation}
{
_a^*D_u^\alpha f(u)= \frac{1}{\Gamma(1-\alpha)}\sum_{k=1}^\infty 
\frac{f^{(k)}(u)}{(k-1)!}(-1)^{k-1} \frac{(u-a)^{k-\alpha}}{k-\alpha}}
\end{equation}
one can determine the value of $u$ that makes $_aD_u^\alpha f(u)=0$. Note that, for a fixed $a$, the interval $[a,t]$ increases with $t$. {An interesting alternative is to consider $a$ as a variable, such as $a = t - L$, where $L$ is a fixed parameter. This results in the function on the interval $[a(t) = t - L, t]$ being considered in the fractional derivative operator.}
Now, memory length is always the same, regardless the value of $t$.
 In that case, the
generalized model was rewritten as
\begin{equation}
u_{k+1}=u_{k}- \omega\hspace{.1cm} [_{u_k-L}D_u^\alpha f(u)](u_k) 
\label{liteq2}
\end{equation}
where $L$ is  memory length.
{Consequently, if $L=h$, $a=u-h$ and $f(u)=(c-u)^2$ thus, the value of $u$, such as $_a^*D_u^\alpha f(u)=0$,} 
is given by 
\begin{equation}
u=c+h(1-\alpha)/(2-\alpha)
\end{equation}
From the above equation, if $h$ goes to zero, then $u=c$ and if $h$ is equal to $u$ then $u=c(2-\alpha)$. The previous algorithms have been widely discussed in the literature, with interesting results \cite{5P}. However, the original question remains: 
the previous algorithms cannot guarantee convergence to the optimal solution \cite{7P,comments}. 
These examples highlight the main disadvantage of the generalized gradient method, { when the fractional operator is substituted for the operator \( \frac{d}{du} \).}

{Having discussed the issue of the extreme point, the next step is to verify the convergence of the methods.} 
In a similar way,
the generalized gradient method can be seen, by the explicit Euler method, as a discretized form of the following ordinary differential equation,
\begin{equation}
\frac{du}{dt}=- \lambda [_{0}^{*}D^\alpha_u f(u)]=F(u)
\label{diferentialgradeq}
\end{equation}
in which $t$ is the iteration time and $\omega=\lambda dt$. 
For instance, consider to minimize $f(u)=(u-c)^2$, using the previous algorithm
and defining $V(u)=(u-c(2-\alpha))^2$, such as 
$V(u) >0$ for all 
$u\neq$ $c(2-\alpha)$ and $V(u)=0$ when $u=c(2-\alpha)$.  
Therefore, $V(u)$ is a positive-definite function. Consequently,  
\begin{equation}
\begin{array}{c}
dV(u)/dt =-2\lambda u(u-c(2-\alpha))[_{0}^{*}D^\alpha_u f(u)]\\
 =-(4\lambda/\Gamma(3-\alpha))(u-c(2-\alpha)^2 u^{1-\alpha}<0  
  \end{array}
\end{equation}
Therefore, for this case, the method has asymptotic convergence to the  
point $u=c(2-\alpha)$, different of the extreme point of the function $f(u)=(u-c)^2$.
In this work, we study a different algorithm to introduce fractional order into GDM and ensure convergence to the extreme point of the function $f(u)$, which will be presented in the following subsection.


\subsubsection{Fractional continuous time algorithm}

The basic idea of this model is then replacing the first-order derivative in time  by a fractional  derivative of 
order $\alpha$, in Equation (\ref{gradeq}), 
so
\begin{equation}
_{0}^{*}D^\alpha_t u=- \lambda\hspace{.1cm} \frac{d}{du}f(u)=F(u)
\label{fracgradeq}
\end{equation}
where $_0^*D_t^\alpha$ is the Caputo operator for the fractional derivative order. 
It is worth highlighting that this proposal is different from other proposals in the literature \cite{4P,5P}. We know few papers that have examined and applied the model that chooses to generalize
gradient method by changing the time derivative instead of gradient \cite{O1,O2}.   

{The right side of the above equation remains as the derivative $df(u)/du$.  
Therefore, the equilibrium point in Equation (\ref{fracgradeq}) is reached when the extremum value of the function $f(u)$ is attained.}
For that reason, we chose to present {Equation (\ref{fracgradeq})} as a generalization of {Equation (\ref{gradeq})}, which works better, as we will see.  Having overcome the issue of the extreme point, the important question that remains is the stability of the equilibrium point.
{For example, consider the minimization of $f(u)=(u-c)^2$ by means of this algorithm,}
\begin{equation}
_{0}^{*}D^\alpha_t u=- \lambda df(u)/du=- \lambda 2(u-c)
\label{fracgradeq1}
\end{equation}
where solution is given by 
\begin{equation}
u(t)=(u(0)-c) E_{\alpha,1}(-2\lambda t^\alpha)+c
\end{equation}
in which $E_{\alpha,1}$ is the Mittag-Leffler function.  In particular, if $0 < \alpha<1$,  we have to specify just one condition, $u(0)$. 
In the above equation, when $t\rightarrow\infty$ then $u(t)\rightarrow c$, which is the extreme point of the function $f(u)$. 
This example shows that,  
using the {Fractional Continuous Time Method},
$u(t)$ converges asymptotically to a  minimal point of the $f(u)$,  whenever 
$\alpha$ is between $0$ and $1$. 
{ Using the Caputo derivative, we can relate the initial condition of the problem to {\it a priori} knowledge of the solution.} 
{ The result obtained by this method was compared with the previous methods in Figure \ref{figure1}.  In this case, we can observe the convergence of the solution to the extreme point of the function $f(u)$. All the cases analysed in Figure 1 use $\lambda=1$.}

In the general case, consider $f(u)$ is  a differentiable 
and $\eta$-strongly convex function, for 
$\eta>0$, in which $u^*$ is the extreme point of $f(u)$. Then $u(t)$, obtained by {Equation (\ref{fracgradeq})}, is unique and converges to $u^*$ at least with the Mittag-Leffler convergence rate.
The Mittag-Leffler stability theorem generalizes Lyapunov theorem 
for fractional order non-linear systems.
 The proof of this result can be found in reference {\cite{O2}}, when $\alpha$ is  restricted in the interval $(0,1]$. Whereas that
$V(u(t))=\frac{1}{2}||u(t)-u^*||^2$,
 using the Inequality (\ref{fracgradin})\cite{ine1} 
\begin{equation}
^*D^\alpha_t V(u(t))\leq -(u^*-u(t)) ^*D^\alpha_t u(t)
\label{fracgradin}
\end{equation}
 and Equation (\ref{fracgradeq}) 
we arrived in
\begin{equation}
{
^*D^\alpha_t V(u(t))\leq (u^*-u(t))\frac{df(u)}{du}}
\end{equation}
Now, using the inequality (\ref{fracgradin2})\cite{ine12} 
\begin{equation}
(u^*-u(t)) \frac{df(u)}{du} \leq f(u^*)- f(u(t)) - \frac{\eta}{2} ||u(t)-u^*||^2
\label{fracgradin2}
\end{equation}
in which $f$ is differentiable and $\eta$-strongly convex, then
we have the result
\begin{equation}
(u^*-u(t)) \frac{df(u)}{du}\leq - \frac{\eta}{2} ||u(t)-u^*||^2 \leq- \eta V(u(t))
\end{equation}
in which $V(u(t))>0$ and $V(u^*)=0$. Finally, using  the Equation (\ref{fracgradeq}) we found that
\begin{equation}
^CD^\alpha_tV(u(t))\leq - \eta V(u(t)) 
\end{equation}
and 
\begin{equation}
V(u(t))\leq V(u(0)){E}_{\alpha,1}(-\eta t^\alpha)
\end{equation}
in which ${E}_{\alpha,1}(-\eta t^\alpha)\rightarrow 0$ when $t\rightarrow \infty$. The Mittag-Leffler stability\cite{Mitag} to solution of the Equation (\ref{fracgradeq}), was demonstrated in literature \cite{7P}, 
with $\alpha$ between $0$ and $1$. 

We do not know of a proof for $\alpha>1$, although simulated results suggest that $u(t)$ converges to $u^*$ when $\alpha>1$, {  how can we seen in {Figure (\ref{figure2a})} at long time. }
To delve deeper into this case, consider to minimize $f(u)=(u-c)^2$ using 
Equation (\ref{fracgradeq1}), with $\alpha$ between 1 and 2,
where solution is given by 
\begin{equation}
u(t)=(u(0)-c) E_{\alpha,1}(-2\lambda t^\alpha)+c
+u^\prime(0)tE_{\alpha,2}(-2\lambda t^\alpha)
\end{equation} 
In this case, two arbitrary initial conditions are necessary, $u(0)$ and $u^\prime(0)$. 
{ Figure (\ref{figure2b}) shows the results found for different values of the $\alpha$ between 1 and 2.} 
{ The performance of the Fractional Continuous Time Method, with \( 1 \leq \alpha \leq 2 \) and different initial conditions, is shown in Figure \ref{figure3}. }
{ For all parameters, the solution exhibits damped oscillations over time. 
For this example, $u(t)$ appears to converge asymptotically as shown in Figure (\ref{figure3}).} The Mittag-Leffler function might have an infinite number of zeros with exception when $\alpha$ is between 0 and 1. We can observe that for each zero of the Mittag-Leffler function, the variable $u$ reaches $u^*$.   

In this study, we propose a global stopping criterion, therefore 
the equations are integrated until a certain objective is achieved, for example, $E<0.01$. 
As we can see in {Figure (\ref{figure2})}, the 
Fractional Continuous Time Method (FCTM)  is capable
of numerically outperforming GDM, { with $\alpha=1.2$ the solution reaches {in $u=3.000$ after} $t=8.39$ arbitrary unit, on the other hand,
with $\alpha=1$ it reaches the solution in $u=2.997$ after $t=32.3$ arbitrary unit.}
FCTM achieved good  precision, $4$ times faster than 
GDM. 
This result is the main interest in the study of the fractional continuous time algorithm. 
However, this depends on the value of $\alpha$.
Considering this example, with $\alpha=0.9$, the performance of FCTM is worse when compared to GDM. 
Hence, how the performance depends on $\alpha$ is a highly important question that is still not well-understood.

In general, FGDM \cite{4P,5P} and Fractional Continuous Time Algorithm \cite{O2} 
have been minimally studied
theoretically and only applied to a 
specific class of functions.  Thus, one main goal of this paper is to solve the practical problems and compare fractional continuous-time algorithm with the usual GDM. In this work, we explore how  the $\alpha$ 
interferes with the convergence rate.
The reference \cite{O2} concludes that $f(\hat{u})$ converges to $f(u^*)$ with at least a rate $O(1/t^\alpha)$, 
in other words, the convergence with $\alpha<1$ is asymptotically slower than that for $\alpha=1$. 
In this context, the variable $\hat{u}$ represents the solution that is being sought.
Nevertheless, for some optimization problems, the observed result is the opposite of this prediction, as we can also see in reference \cite{O2}.
Despite some theoretical results for $0<\alpha<1$, we are not aware of theoretical studies considering $\alpha>1$. In this work, we show some numerical results that can be used to 
suggest directions for future theoretical research.
{In Figure (\ref{figure4}), we} can see that energy function $V(t)=(u(t)-3)^2$, for $\alpha>1$, does 
not decay as expected for $\alpha$ between $0$ and $1$.

The Hopfield neural network (HNN) is one of the most used neural network architectures, it has been used to solve ill-posed problems with great success. Some works in recent years have incorporated fractional-order derivative with respect to time into original Hopfield neural network model.  In reference Tavares et al. (2022)\cite{Nel3} fractional-order derivative with respect to time was included on Hopfield neural network equations, obtained from Lyapunov function defined by 2-norm of the residual function. In the HNN model, the state of the neurons $j$ in time $t$, $u_j(t)$, is a function of the $g_j(t)$,
such that $u_j(t)=h^{-1}(g_j(t))$ and $g_j(t)=h(u_j(t))$, with $j=1...m$ and $m=dim({ g})$. In this case, if  $g_j(t)=u_j(t)$, with $j=1...m$, 
we recover the Equation (\ref{fracgradeq}). The question about how fractional order has an effect on the solution found was discussed using the Mittag-Leffler
stability of the fractional-order Hopfield neural network as criterion. This is analogous to the one presented here. The result found by the 
Fractional-order Hopfield neural network model has been achieved in less time when compared to the Hopfield neural network model.
\begin{center}
\begin{figure}[h]
\centering
\includegraphics[scale=.6]{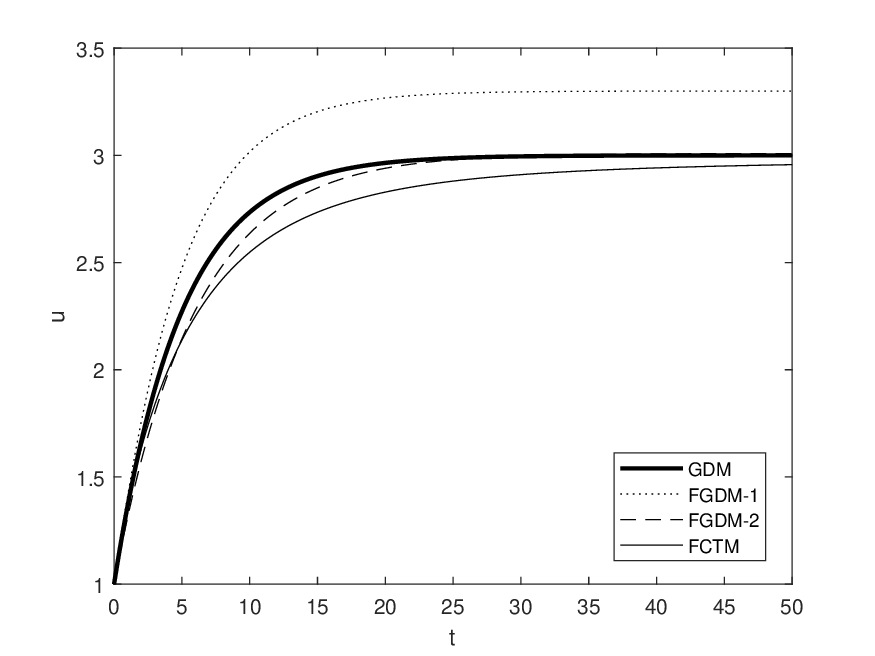}
\caption{The optimal value obtained at each iteration, using Gradient Descent Method (FGDM) (thick solid line). The solutions obtained by GDM-1 \cite{4P} and (FGDM-2) \cite{5P}, both  
 with $\alpha=0.9$, is represented by the dotted line and dashed line, respectively.
The   thin continuous line shows the result found by Fractional Continuous Time Method (FCTM), 
with $\alpha=0.9$. The expected result is 3 ($c=3$).}
\label{figure1}
\end{figure}
\end{center}

\begin{figure}[htbp]
     \centering
     \subfigure[Result obtained with $\alpha=1$ (thick solid line), 
$\alpha=0.9$ (dotted line), $\alpha=0.7$ (dashed line) and
$\alpha=0.5$ (thin continuous line). \label{figure2a}]
{
\includegraphics[width=.4\textwidth]{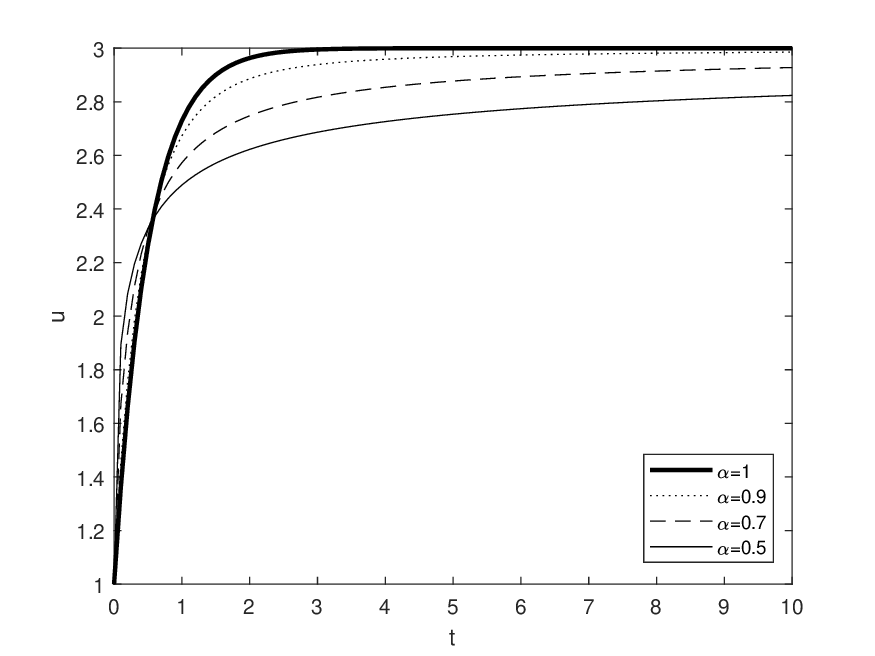}
}
\quad
      \subfigure[Result obtained with $\alpha=1$ (thick solid line), 
$\alpha=1.2$ (dotted line), $\alpha=1.5$ (dashed line) and
$\alpha=1.7$ (thin continuous line). \label{figure2b}]
{
\includegraphics[width=.4\textwidth]{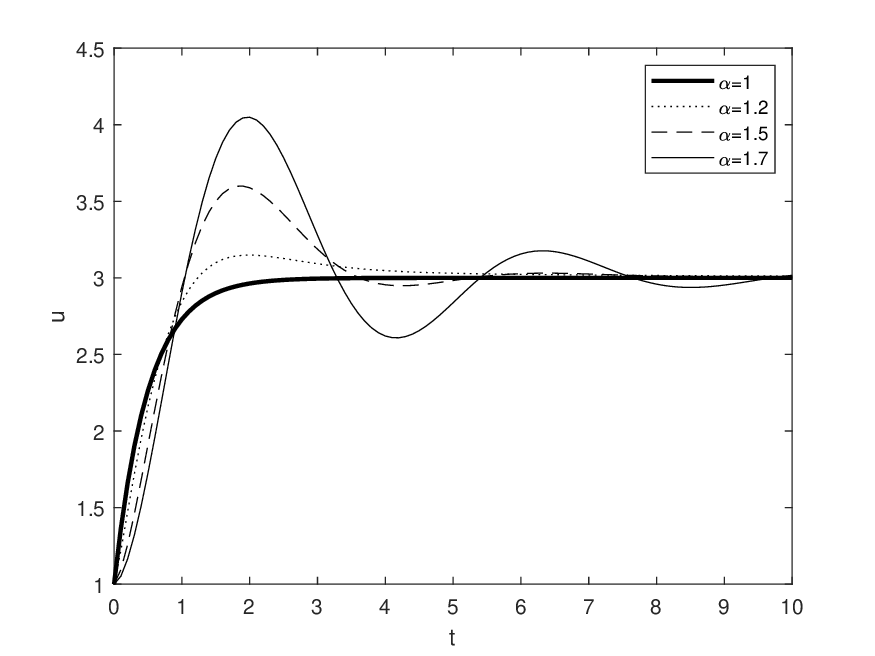}
}
 \caption{ The solutions obtained by  Fractional Continuous Time Method (FCTM), 
with different values of fractional order $\alpha$, for the problem of minimizing the function $f(u)=(u-3)^2$. }
\label{figure2}
\end{figure}

\clearpage
\begin{center}
\begin{figure}[h]
\centering
\includegraphics[scale=.6]{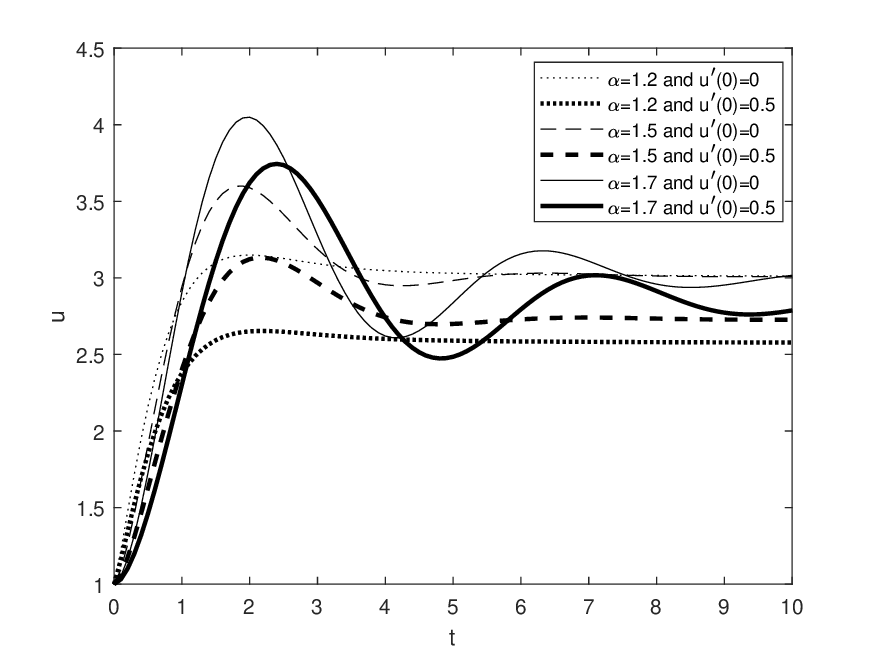}
\caption{The performance of the Fractional Continuous Time Method, with $1\leq\alpha\leq 2$ and different initial conditions: for all $u(0)=1$ and $u^\prime(0)=0$ (thin line) or $0.5$ (thick line).
Result obtained with $\alpha=1.2$ (dotted line), $\alpha=1.5$ (dashed line), 
and $\alpha=1.7$ (continuous line).
}
\label{figure3}
\end{figure}
\end{center}

\begin{figure}[htbp]
     \centering
     \subfigure[Result obtained with $\alpha=0.9$ (thick doted line), $\alpha=1$ (thick solid line), 
$\alpha=1.2$ (dotted line), $\alpha=1.5$ (dashed line) and
$\alpha=1.7$ (thin continuous line). \label{figure4a}]
{
\includegraphics[width=.4\textwidth]{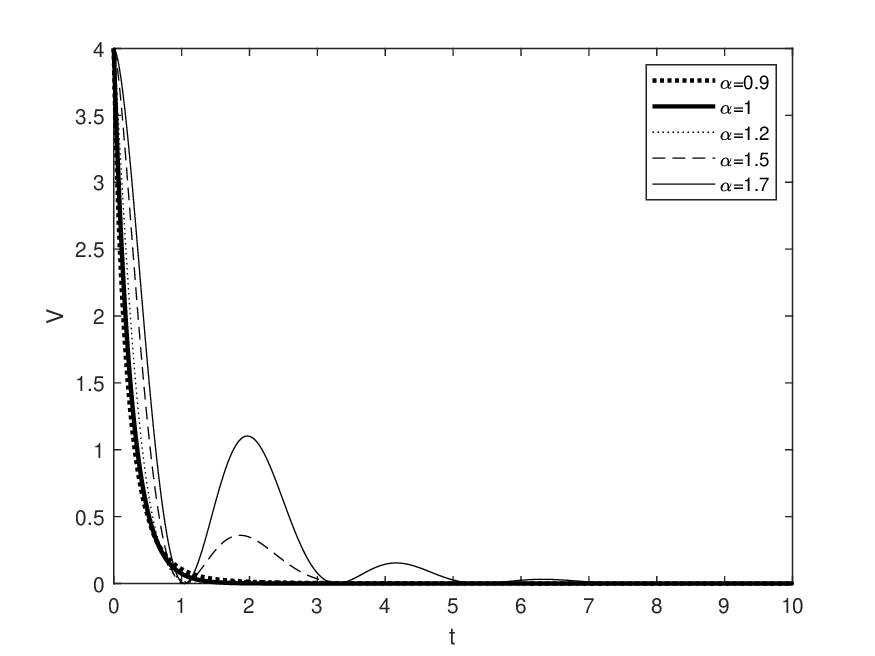}
}
\quad
      \subfigure[Zoom of Figure (\ref{figure4a}) for $0 \leq V \leq 1.4$ and $0 \leq t \leq 5$.  \label{figure4b}]
{
\includegraphics[width=.4\textwidth]{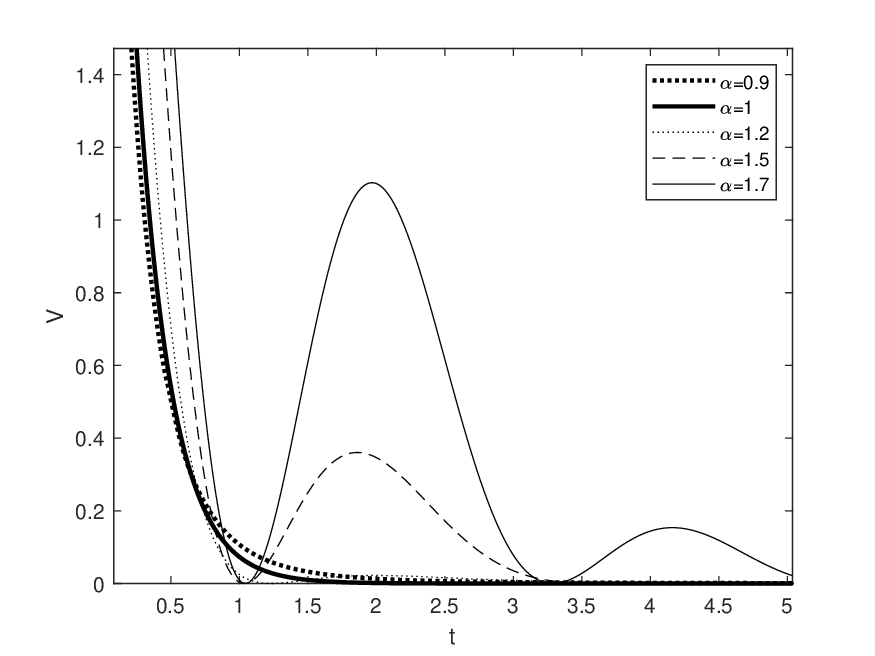}
}
 \caption{The error $V(u)=(u(t)-u^*)^2$ as a function of time. This result is obtained by  Fractional Continuous Time Method (FCTM), for the problem of minimizing the function $f(u)=(u-3)^2$. }
\label{figure4}
\end{figure}

\section{Applications}

\subsection{Interpolation polynomial of degree $m$}

The next example discusses the case of the Vandermonde matrix, which is a well-known ill-posed problem.
Given the values of a function $g(x)$ for two different values of $x$, for example $x_0$ and $x_1$, we can approximate $g(x)$ by a polynomial function
of degree 1,
$p_1(x)=u_0k+u_1$, which satisfies these conditions: $p_1(x_0)=g(x_0)$ and $p_1(x_1)=g(x_1)$.
In this case, the desired  solution  (values of $u_0$ and $u_1$) is obtained by solving the system: $p_1(x_0)=u_0x_0+u_1=g(x_0)$ and $p_1(x_1)=u_0x_1+u_1=g(x_1)$.

The interpolation polynomial of degree $m$ can be written as $P_m=u_m+u_{m-1}x_1+u_{m-2}x_2+...+u_0x_m$. In order for $P_m(x)$ to replace $g(x)$, the coefficients $u_j$ must be determined for all $1\leq j\leq m$. To fit these coefficients to the data, $P_m(x_j)=g(x_j)$ must be satisfied for all $j$ between $1$ and $m$. In this case, the coefficients
are  obtained by solving the linear system
${\bf X}{\bf u}={\bf g}$, i.e.
\begin{equation}
\left[
\begin{array}{ccccc}
x_0^m &...& x_0^2 &x_0  & 1\\
x_1^m &...& x_1^2 &x_1  & 1\\
...& ... &... &... &...\\
x_m^m &...& x_m^2 &x_m  & 1\\
\end{array}\right]
\left[\begin{array}{c}
u_0\\
u_{1}\\
...\\
u_m
\end{array}\right]=
\left[\begin{array}{c}
g_0\\
g_1\\
...\\
g_m
\end{array}\right]
\end{equation}
The above problem can be seen as an optimization problem because there is an evaluation 
function $f({\bf u})=||{\bf Xu}-{\bf g}||^2$ which gives a score to the candidate function ${ u}$.
The majority of recent applications are limited to optimization problems with a small number of variables, typically $n=1$ or $2$.
To test the performance of the { FCTM}, 
the Vandemonde matrix 11 $\times$ 11, where  $0<x<1$, was used. 
The predictor-corrector method of Adams-Bashforth-Moulton described in \cite{PECE1,PECE2,PECE3,PECE4,PECE5} (PECE method) and implemented 
in \cite{fde12} (fde12 implementation) was used to 
solve Equation (\ref{fracgradeq}), with $\lambda=0.001$. The stability properties of the method implemented have been studied in \cite{PECE1,PECE2,PECE3,PECE4,PECE5}.

The comparison between $u_j$, for all $1\leq j\leq m$,
over time, obtained using $\alpha=1$ and $\alpha=1.2$,  is shown in Figure (\ref{figure5a}).
The results were obtained from an initial condition of ${ u}(0)=0$ and ${ u}^\prime(0)=0$ (when $\alpha >1$).
Figures (\ref{figure5b}), (\ref{figure5c}) and (\ref{figure5d}) show the residual norm for FCTM with $\alpha=0.8$, $1.2$ and 
 $1.4$.  The residual norm, when FCTM is used with $\alpha=0.8$, decreases more slowly than that obtained by  GDM (or FCTM with $\alpha=1$) at the same interation time (arbitrary units).
On the other hand, when FCTM is used with $\alpha=1.2$ or $1.8$, the residual norm decreases faster than the value obtained by GDM. Using FCTM, with  $\alpha=1$ (which is the same as using GDM), $||{\bf Xu}-{\bf g}||_{t=50000}= 1.7 \times 10^{-5}$  about 94 times greater than residual norm, if $\alpha=1.2$ is used, $||{\bf Xu}-{\bf g}||_{t=50000}= 1.8 \times 10^{-8}$.  { These results are presented in Table \ref{tab1}}. 
{ The performance of the Fractional Continuous Time Method, with \( \alpha =1.2 \) and different initial conditions, is shown in Figure \ref{figure6}. }

{ It is important to note that the numerical approach for integer-order differential equations in GDM was the Runge-Kutta method with a variable time step (ode45 implementation) \cite{1ode,2ode}. For a system of fractional differential equations, another method is required, as discussed earlier. Fractional-order systems are non-local, and the corresponding numerical methods lead to very complex schemes, where every computed step relies on all previously computed steps.
Therefore, solving a fractional-order system numerically over a large time interval incurs significantly higher computational costs compared to classical ordinary differential equations

In order to clarify the question, the computational cost was computed for a processing time of \( t = 5000 \) a.u. (arbitrary unit). The processing time using the Runge-Kutta method (\(\approx\) 0.75 s) is approximately 20 times lower than the time using the numerical approach for fractional-order differential equations with \( \alpha = 1 \) (\(\approx\) 15 s). A similar time is observed for other values of \( \alpha \).

The computational cost was calculated using an Intel(R) Core(TM) i5-2400 CPU @ 3.10 GHz. Even so, the fractional algorithm offers a better cost-benefit ratio. The residual decreases by almost 94 times for \( \alpha = 1.2 \) and about 629 times for \( \alpha = 1.4 \), while the computation time is only 20 times greater. These results are presented in Table 1. Therefore, if the fractional order is chosen appropriately, the {FCTM} model achieves the desired solution in less computing time than the GDM model.} 

\begin{figure}[htbp]
     \centering
     \subfigure[Enhancement of the solution ${ u}$ with respect to time, using FCTM {\bf with 
      $\alpha=1.2$ (dotted line) and $\alpha=1$ (solid line)}. \label{figure5a}]
{
\includegraphics[width=.4\textwidth]{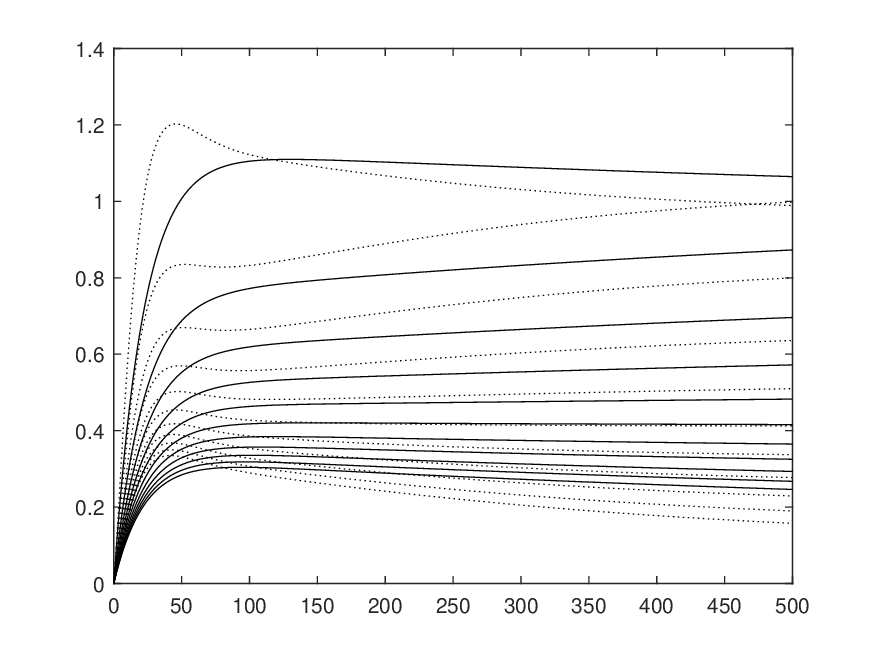}
}
\quad
      \subfigure[Residual $||{\bf Xu}-{\bf g}||$ through time, using FCTM, for a linear system $11\times 11$ and fractional order {equal to $0.8$ (dotted line) and $\alpha=1$ (solid line)}. \label{figure5b}]
{
\includegraphics[width=.4\textwidth]{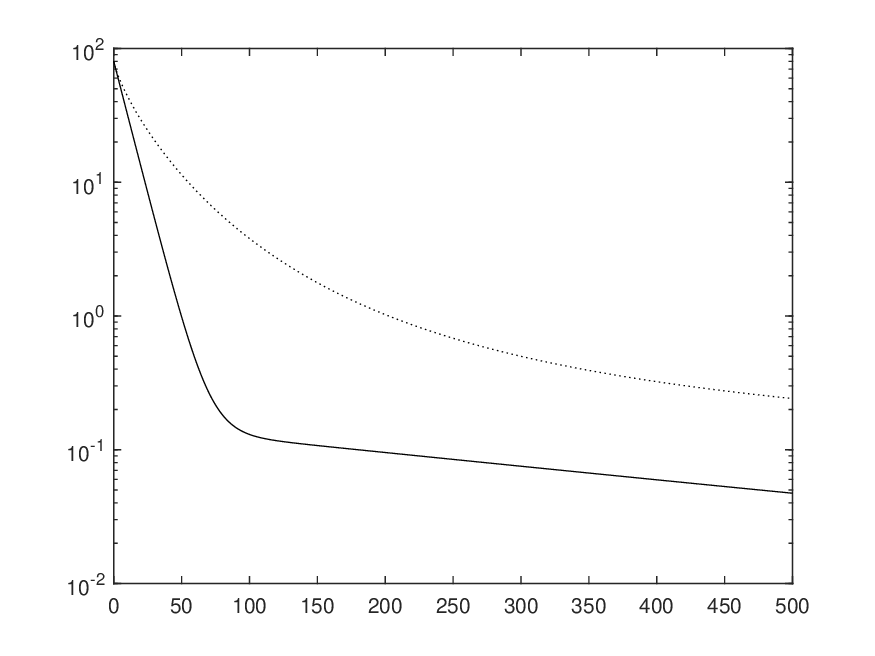}
}
\quad
      \subfigure[The residual function $||{\bf Xu}-{\bf g}||(t)$ obtained by  Fractional Continuous Time Method (FCTM) {with 
      $\alpha=1.2$ (dotted line) and integer order model is represented by solid line}. \label{figure5c}]
{
\includegraphics[width=.4\textwidth]{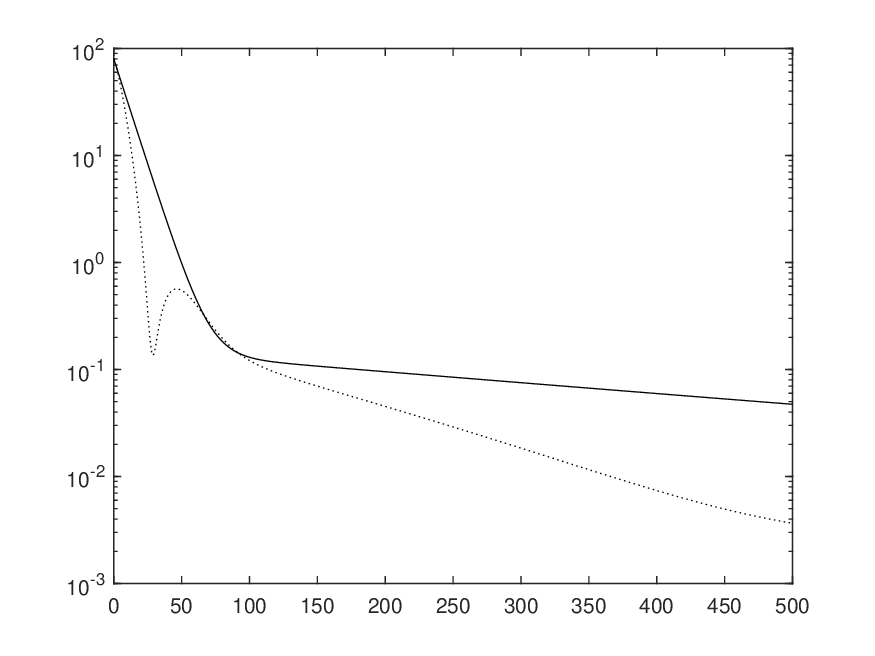}
}
\quad
      \subfigure[Result ${\bf u}(t)$ {with $\alpha=1.4$ is represented by dotted line and 
integer order model is represented by solid line.} \label{figure5d}]
{
\includegraphics[width=.4\textwidth]{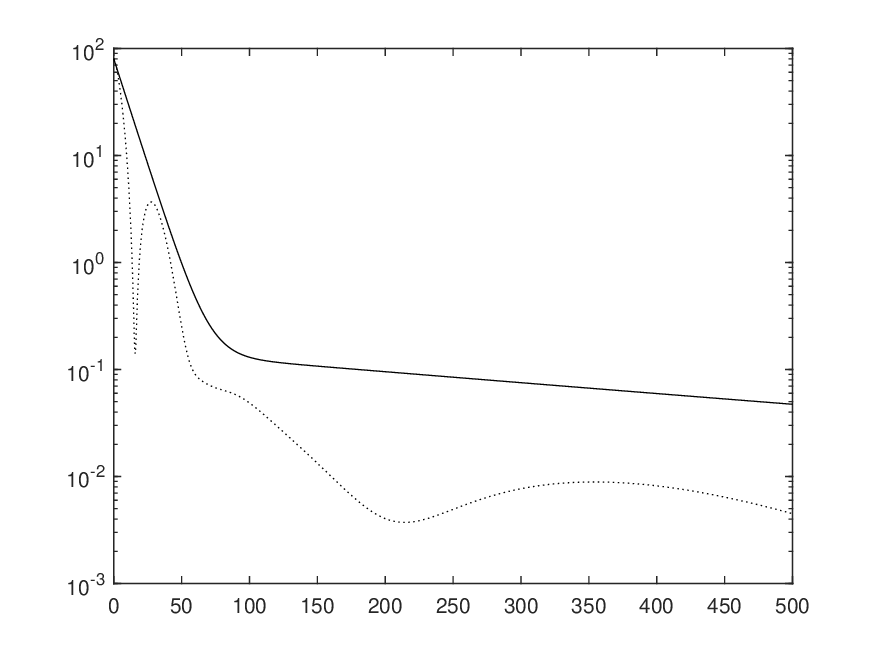}
}
 \caption{The performance of the Fractional Continuous Time Method  to solve the interpolation polynomial problem.}
\label{figure5}
\end{figure}

\begin{figure}[htbp]
     \centering
     \subfigure[Enhancement of the solution ${\bf u}$ with respect to time, using FCTM with 
      $\alpha=1.2$ and initial condition ${\bf u}_0={ 1}$ {(dotted line). 
The integer orde model is represented by solid line.}      
      \label{figure5ax}]
{
\includegraphics[width=.4\textwidth]{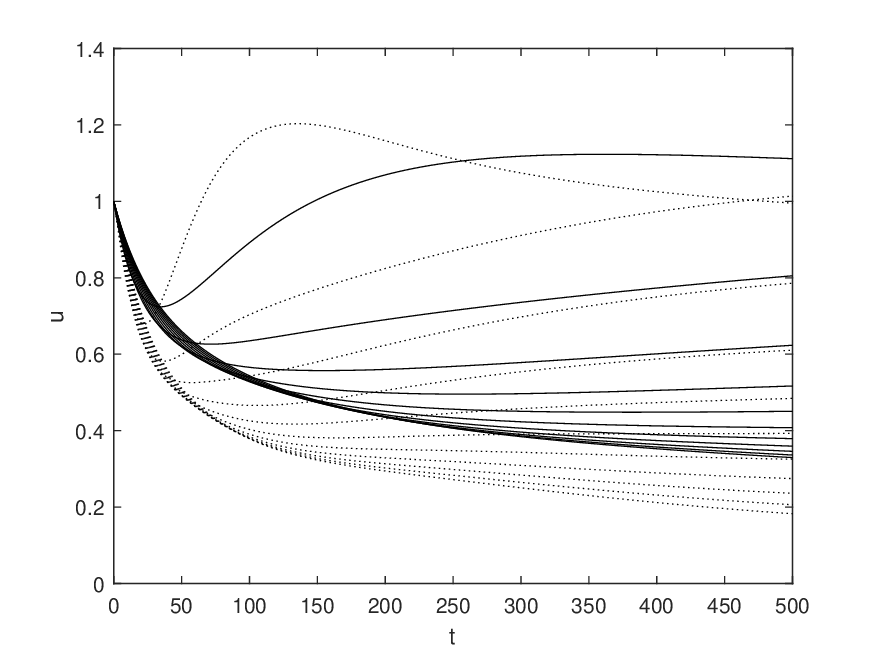}
}
\quad
      \subfigure[Residual $||{\bf Xu}-{\bf g}||$ through time, using FCTM with 
      $\alpha=1.2$ and initial condition ${\bf u}_0={ 1}$ {(dotted line). The integer orde model is represented by solid line.} \label{figure5bx}]
{
\includegraphics[width=.4\textwidth]{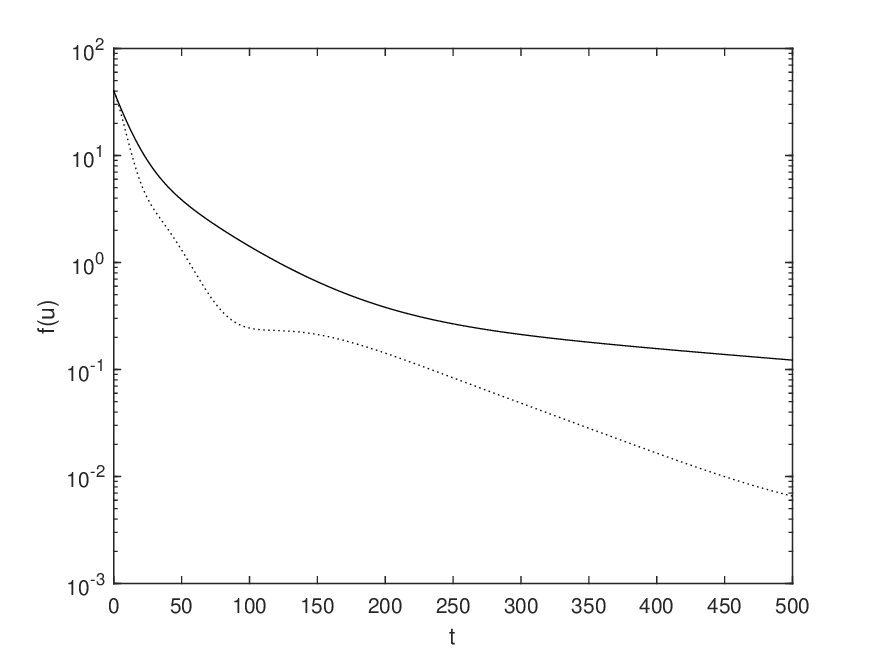}
}
\quad
      \subfigure[Enhancement of the solution ${\bf u}$ with respect to time, using FCTM with 
      $\alpha=1.2$ and random initial condition {is represented by dotted line. 
 The solid line shows the result obtained for the integer-order model.}      
       \label{figure5cx}]
{
\includegraphics[width=.4\textwidth]{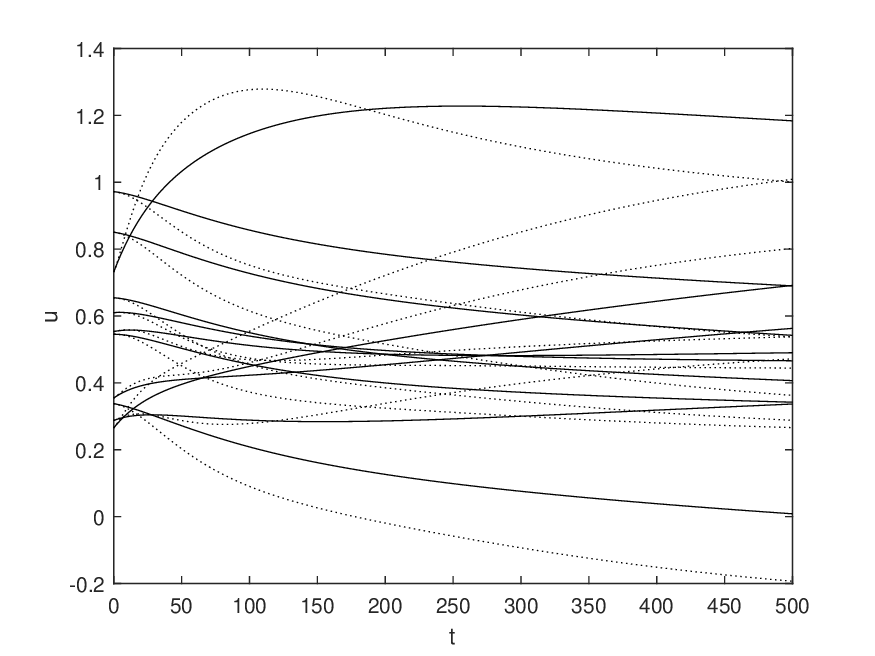}
}
\quad
      \subfigure[Residual $||{\bf Xu}-{\bf g}||$ through time, using FCTM with 
      $\alpha=1.2$ and random initial condition {(dotted line).
The solid line shows the result obtained for the integer-order model.}      
        \label{figure5dx}]
{
\includegraphics[width=.4\textwidth]{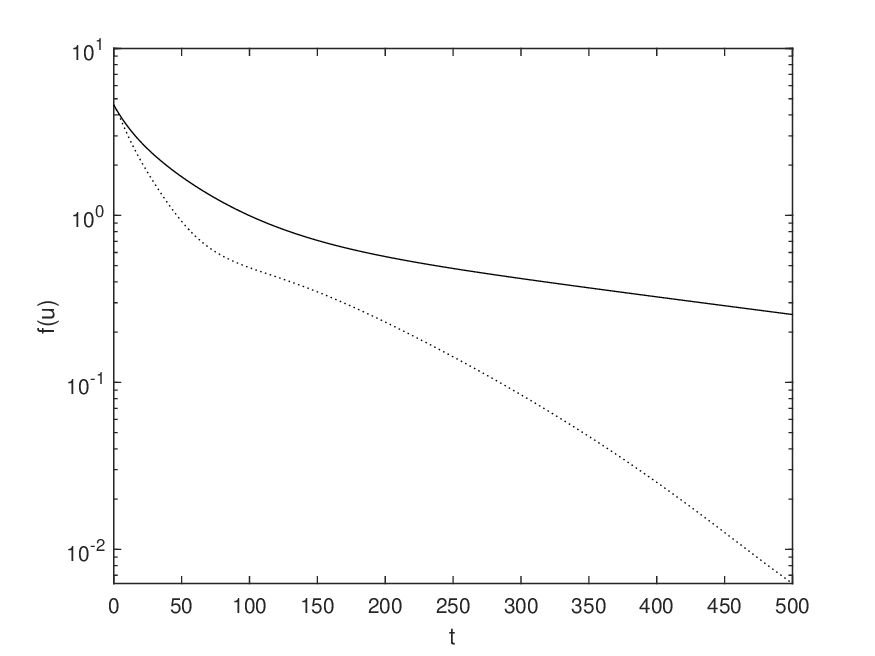}
}
 \caption{The performance of the Fractional Continuous Time Method  to solve the interpolation polynomial problem.}
\label{figure6}
\end{figure}

\begin{table}[]
\centering
\caption{
Comparison of the results using GDM model ($\alpha=1$) and FCTM model.}
\label{tab1}
\begin{tabular}{c|c|c|c|c|c}
\hline\hline
  $\alpha$   & $t^\alpha_{<0.1}$, a.u.   &  $t_{<0.01}^{\alpha}$, a.u. &  $t_{<0.001}^{\alpha}$, a.u. &
  $||{\bf Xu}-{\bf g}||^\alpha_{t=50000}$ & $\left(\frac{||{\bf Xu}-{\bf g}||^1}{||{\bf Xu}-{\bf g}||^\alpha}\right)_{t=50000}$\\ \hline\hline
 0.8 & 1139 & 9444& $>50000$  & 1.5$\times 10^{-3}$ & 0.011 \\ \hline
 1.0 & 180  & 1245& 6877  & 1.7$\times 10^{-5}$ & 1 \\ \hline
 1.2 & 113  & 366 & 1758  & 1.8$\times 10^{-7}$ & 94 \\ \hline
 1.4 & 58   & 160 & 670   & 2.7$\times 10^{-8}$ & 629 \\ \hline
 1.6 & 35   & 84  & 325   & 3.3$\times 10^{-8}$ $^*$ & 515 \\ \hline\hline
\end{tabular}\\
\hspace{-6cm}$^*$ {\small The result begins to show numerical instability}
\end{table}

\subsection{Thomson problem}

The problem of distributing point charges on a sphere has a long history, starting with Thomson in 1904 when he proposed his atomic model \cite{1T}. The objective of the Thomson problem is to determine the minimum electrostatic potential energy configuration of 
$N$ charged particles, all with equal charges, constrained to the surface of a sphere with the radius equal to 1. The electrostatic interaction energy between each pair of charged particles is given by Coulomb's law, and the total electrostatic potential energy 
may then be expressed as the sum of all pair interaction energies, such that 
\begin{equation}
f({\bf u})=\sum_i^N \sum_{j=i+1}^N r_{ij}^{-1}
\end{equation}
in which $r_{ij}=[(x_i-x_j)^2+(y_i-y_j)^2+(z_i-z_j)^2]^{-1/2}$, $x_i=\sin(\phi_i)\cos(\theta_i)$, 
$y_i=\sin(\phi_i)\sin(\theta_i)$ and $z_i=\cos(\phi_i)$. The  
 discrete charges constrained
to a spherical volume of space, such that $x_i^2+y_i^2+z_i^2=1$. This simplifies the calculations by reducing the search space  from 3N dimensions down to 2N dimensions,  
in which ${ u}=[\theta\hspace{.1cm} \phi]$.
The global minimization of total electrostatic potential energy
over all possible configurations  can be found by 
numerical minimization algorithms, such as GDM. 

The solutions of the Thomson problem for $N = 4$, $5$, $6$ and $12$ charged particles are well-known to mathematicians and chemists \cite{2T}. The geometric solution of the Thomson problem for $N = 4$, $6$, and $12$ electrons is a Platonic solid in which the faces are all congruent equilateral triangles.  Apart this chemical interest, the Thomson problem has been used as a benchmark problem, with large values of $N$, in testing global optimization algorithms \cite{3T}.  
The global minima for the Thomson Problem with small values of $N$ can be found, for instance, in references  \cite{4T,5T}. The lowest found minima for the Thomson problem for some larger $N$ values  are available in references \cite{6T,7T}. In this section, we present numerical solution of  Thomson problem {with $N=4$, $5$, $6$ and $12$, using GDM and FCTM.} { The Thomson problem with 
$N=12$ has $24$ parameters to be optimised. }

{The Figure} (\ref{figure7a}) shows that $f({ u})$, when FCTM is used with $\alpha=0.7$, initially 
drops to less
 than that obtained with GDM (or FCTM with $\alpha=1$).  Then there is a change in the 
  in the decay behaviour, and $f({ u})$ falls more slowly than with $\alpha=1$.
  In this study, we propose a global stopping criterion, therefore 
the equations are integrated until a certain objective is achieved. 
As we can see in {Figure (\ref{figure7b})}, the 
FCTM is capable
of numerically outperforming GDM, with $\alpha=0.7$, for a small number of iterations. However, the Fractional continuous time Method (FCTM) achieved good  precision faster than 
Gradient Descent Method (GDM). This result can be used as an initial condition in more efficient methods. 

{ The results presented for the computation time show, as expected, that the use of the ode45 routine is more efficient than the fde12 routine (with \( \alpha = 1 \)). It can also be observed that the computational time ($\tau$) is approximately the same for any \( \alpha \). All the results in Table \ref{tab2} are computed for an iteration time ($t$) of 1500 a.u. and a step size of \( h = 1 \times 10^{-5} \). For the first time in the literature, we present a comparison of the computation time between the classical gradient method and its generalisations to fractional order. In Table \ref{tab2}, for \( N = 12 \), we observe a gain in the cost function of \( \frac{49.165253058 - 49.207587517}{49.165253058 - 49.167200926} = 21.7 \) times, while the computation time increases by \( \frac{8.2}{0.61} = 13.4 \) times. Therefore, there is a real gain in using the FCTM. Figure \ref{figure7b} presents the optimal geometry for the case $N = 12$, where the solution corresponds to a regular icosahedron.}

\begin{center}
\begin{figure}[h]
\centering
\includegraphics[scale=.6]{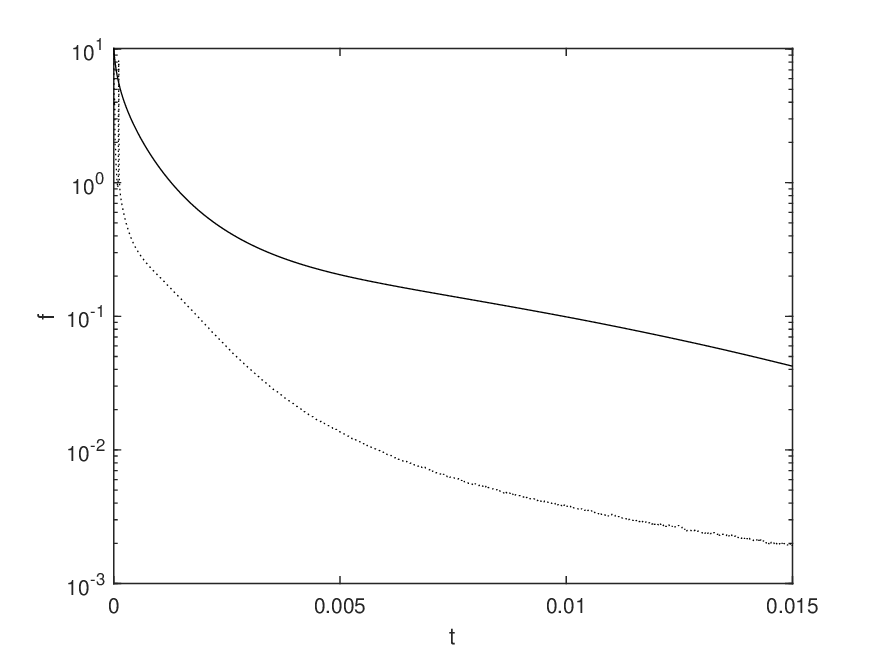}
\caption{The residual over time, using FCTM and a fractional order of \( \alpha = 0.7 \), is shown by dotted line. The result obtained with \( \alpha = 1 \) is represented by the solid line. }
\label{figure7a}
\end{figure}
\end{center}

\begin{center}
\begin{figure}[h]
\centering
\includegraphics[scale=.6]{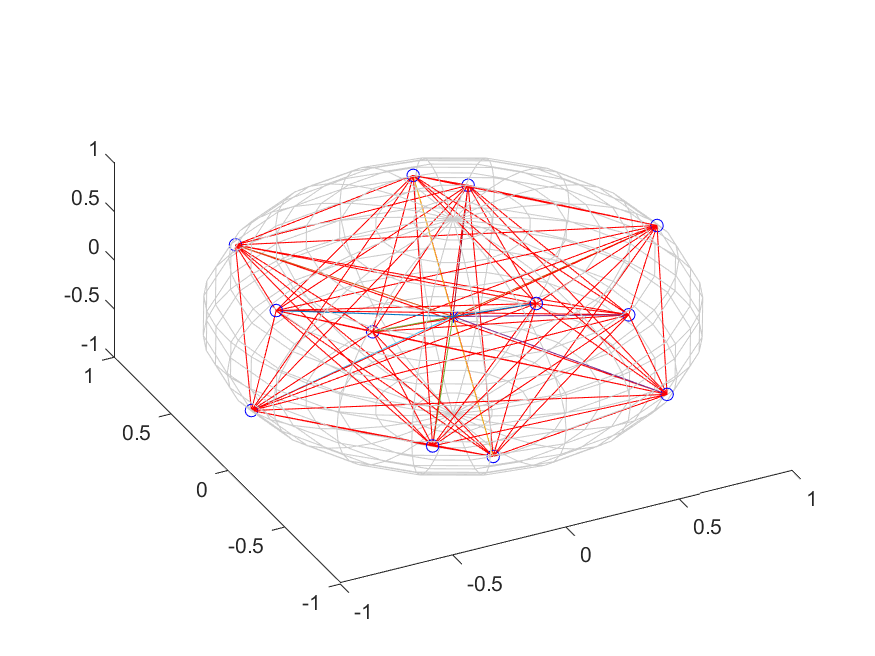}
\caption{Minimum electrostatic potential energy configuration of 
$6$ charged particles, all with equal charges, constrained to the surface of a sphere with the radius equal to 1. The geometry found is that of a regular icosahedron.}
\label{figure7b}
\end{figure}
\end{center}

\begin{table}[]
\centering
\caption{
Comparison of the results using GDM model ($\alpha=1$) and FCTM model ($\alpha=0.7$) for the Thomson Problem.\vspace{.2cm}}
\label{tab2}
\begin{tabular}{c|c|c|c|c|c|c}
\hline\hline
$N$  &    $E_\text{reference}$  &  $E_1$     & $\tau_1$,s(ode45) & $\tau_1$,s(fde12)  & {$E_{\alpha=0.7}$}   & {$\tau_{\alpha=0.7}$,s}
\\ \hline\hline
4  &  3.674234614  & 3.733304005  & 0.20 & 1.4  &  3.678389136  & 1.6\\
5  &  6.474691495  & 6.522966573  & 0.16 & 1.9  &  6.476946521  & 2.0\\
6  &  9.985281374  &10.096028159  & 0.18 & 2.4  &  9.992883145  & 2.5\\
12 & 49.165253058  &49.207587517  & 0.61 & 8.1  & 49.167200926  & 8.2\\ \hline\hline
\end{tabular}\\
\end{table}

\section{Conclusions}

{Many problems in science and engineering involve finding a set of parameters $x$ that minimizes an objective function $f(x)$. The 
gradient descent  method is one of the simplest and most commonly used methods to solve problems posed in this form.  In general, the gradient descent method is slow to converge close to the extreme point. In this perspective, fractional
calculus has been explored as a promising tool to improve the ordinary gradient descent method. Meanwhile, using the fractional order of the derivative can have unexpected consequences compared to the integer order of the derivative. 

{ This article compares the different methods presented in the literature with regard to the convergence of the method, convergence to the extreme point and convergence rate. The result is that, in general, methods that generalize the gradient to fractional order 
are inefficient in obtaining an efficient optimization method.
As discussed here  the point at which fractional gradient descent converges is highly dependent on
the choice of derivative definition, integral limit and fractional-order. Therefore, convergence to the extreme point of the function that you want to minimise cannot be guaranteed. This point is highlighted with several examples.}

 To avoid these difficulties, we have chosen the Fractional Continuous Time algorithm to generalise the gradient method. Changing the time derivative instead of the gradient guarantees convergence to the extreme point of the coast function. This route has been little explored in the literature and, as presented here, has more promising characteristics. In this case the issue of the extreme point is overcome, meanwhile remains the issue about the stability of the equilibrium point.{  The Mittag-Leffler stability to solution obtained by FCTM is demonstrated to $\alpha$ between 0 and 1. In this work, computational simulation suggest that  optimization parameter, obtained by fractional continuous time method, converges to extreme point of coast function, when  fractional-order  is between 1 and 2.} 
{ There is still no rigorous proof for this result in the literature.}

Moreover, the examples shown in this paper illustrate that, for some fractional order, the convergence speed of 
fractional continuous algorithm is faster than gradient descent method, if the fractional order is chosen appropriately.}
{ The first example
deals with the linear system and second with non-linear problem.  In two of these examples, the rate of convergence was
greater with $\alpha \neq 1$. }
{ The key questions that were addressed in this work are important for the use of FCTM model in context
of the experimental Physical-Chemistry field. In general, previous studies are restricted to mathematical
questions and examples.} 
 Despite the positive results, future work should be carried out to clarify why fractional gradient descent improves the usual method and how to choose the fractional order.}

\section*{Acknowledgments}

This study was financed in part by the Coordena\c c\~ao de
Aperfei\c coamento de Pessoal de N\'ivel Superior (CAPES) - Finance Code 001, and Funda\c c\~ao de Amparo \`a Pesquisa do Estado de Minas Gerais (FAPEMIG). We would like to thank the comments of  
Jos\'{e} Claudinei Ferreira, Institute of Exact Science, Universidade Federal de Alfenas (UNIFAL), 
which helped to improve this paper. {We sincerely thank the anonymous reviewers for their insightful comments and suggestions that helped improve this work.}











\end{document}